\documentclass[12pt,a4paper]{article}
\pdfoutput=1
\usepackage[italian]{babel}
\usepackage{amsmath}
\usepackage{amssymb}
\usepackage{latexsym}

\numberwithin{equation}{section}

\newtheorem{definition}{Definizione}[section]
\newtheorem{corollary}{Corollario}[section]
\newtheorem{proposition}{Proposizione}[section]
\newtheorem{lemma}{Lemma}[section]

\def\Q{\mathbb Q}
\def\N{\mathbb N^*}

\def\R{\mathbb R}
\def\CB{\mathcal B}
\def\CA{\mathcal A}

\def\CC{\mathcal C}

\def\ag{$\rm{\grave{a}}$}

\def\ea{$\rm{\acute{e}}$}
\def\a{\`a\ }
\def\e{\`e\ }          
\def\io{\`{\i}\ }
\def\o{\`o\ }
\def\u{\`u\ } 
\def\E{\`E\ }  
\def\df{\displaystyle\frac}
\def\bx{\hspace{\stretch{1}} $\Box$ \vspace{1mm}}

\def\it{\textit }
\def\rm{\mathrm } 

\begin{document}
\begin{center}
\Large
\textbf{Una costruzione del sistema dei numeri \vspace{6mm} reali} \\
\normalsize
\textbf{Maria Rosaria\; ENEA,\; \ Donato\; \vspace{8mm} SAELI}
\begin{minipage}{120mm}
\small
\textsc{Abstract.} There are many ways to construct the field $\, \R \,$ of real numbers. The most important and famous of these employ Cauchy sequences (Cantor) or cuts (Dedekind) in the field $\, \Q \,$ of rational numbers. These constructions sometimes overlook important details and often are complicated and cumbersome. In this note, the authors propose an essential, clear and rigorous construction of $\, \R \,$ from the stucture $\, \Q_+^{^*} \,$ of positive rational numbers by the key notion of initial \vspace{2mm} segment.  \vspace{2mm} segments. \\
\textsc{MSC:} \vspace{6mm} 97F50
\end{minipage}
\end{center}
\section{\large Introduzione}
\normalsize
Nella seconda met\a dell'Ottocento si avverte pi\u che mai la mancanza di un fondamento rigoroso del calcolo infinitesimale: l'uso di argomenti geometrici, anche se didatticamente utili, non potevano essere considerati un fondamento \textit{scientificamente} valido per il calcolo.
Da qui la necessit\a di \textit{scoprire negli elementi dell'aritmetica la vera origine} del calcolo.\\
Contributi decisivi all'\textit{aritmetizzazione dell'analisi} furono dati da Weierstrass, Dedekind e Cantor; contributi che affondano le loro radici in mezzo secolo di ricerche sulla natura delle funzioni e dei numeri. \\
Dedekind aveva rivolto la sua attenzione ai numeri irrazionali sin dal 1858: i numeri razionali non consentivano di descrivere aritmeticamente le propriet\a della retta, occorreva creare nuovi numeri che consentissero di ottenere un campo \textit{altrettanto continuo come lo \e la retta}. Una propriet\a caratteristica della continuit\a, che servisse da base per costruire il nuovo campo, viene da lui formulata nel seguente principio:\\
"\textit{Se una ripartizione di tutti i punti della retta in due classi \e di tale natura che ogni punto di una delle due classi stia a sinistra di ogni punto dell'altra allora esiste uno e un solo punto dal quale questa ripartizione di tutti i punti in due classi, o questa decomposizione in due parti, \e prodotta}" \\
Ogni numero razionale produce una tale \textit{sezione}, cos\io come le chiama Dedekind, ma esistono anche infinite sezioni non prodotte da numeri razionali; l'esempio pi\u semplice \e la sezione $(A, B)$ definita ponendo in $A$ tutti i razionali negativi, lo zero e tutti i razionali positivi il cui quadrato sia minore di 2 e nella classe $B$ tutti i rimanenti razionali.\\
"\textit{Orbene}- scrive Dedekind -\textit{ogni volta che \e data una sezione che non sia prodotta da un numero razionale noi creiamo un nuovo numero, un numero irrazionale $\alpha$ che noi consideriamo completamente definito da questa sezione; noi diremo che il numero $\alpha$ corrisponde a questa sezione e che esso la produce.}"\\
Dedekind dimostra le propriet\a di ordinamento del nuovo sistema di numeri, la continuit\a e ne definisce le usuali operazioni. \\
Contemporaneamente Cantor, che aveva seguito i corsi di Weierstrass a Berlino, studiando il problema dell'unicit\a della rappresentazione in serie di Fourier di una funzione, formula un assioma che Dedekind ritiene essere equivalente al suo assioma di continuit\a. \\
Cantor considerava successioni infinite di numeri razionali $\{a_n\}$ soddisfacenti alla seguente condizione: per ogni $\epsilon>0$, esiste un intero $\bar{n}$ tale che per $n\geq \bar{n}$ e $m$ intero qualunque, sia $|a_{n+m}-a_n|<\epsilon$. Ad ogni successione che soddisfa tale condizione, detta successione \textit{fondamentale} (o "di Cauchy"), rimane \textit{associato} un numero\footnote{Cantor scrive "\textit{ha un limite determinato} $b$".} $b$: per questi numeri Cantor definiva il concetto di uguaglianza e mostrava come si potessero estendere ad essi le usuali operazioni.\\
Fissata un'unit\a di misura, ogni punto della retta poteva essere associato ad uno di questi numeri; il viceversa era assicurato dal seguente assioma,
"\textit{ad ogni numero corrisponde un ben determinato punto della retta, la cui coordinata \e uguale a quel numero}".\\
Ancora oggi questi due metodi di costruzione dei reali sono i pi\u usati nell'insegnamento. In verit\a tale costruzione viene spesso solo abbozzata a causa delle difficolt\a tecniche che si incontrano nel definire le operazioni e nel dimostrare le propriet\a di cui godono.\\ 
Nel testo di algebra di Seth Warner, viene presentata, nel capitolo \textit{The real and complex number fields}, la costruzione dei reali secondo Cantor, ma alla fine del paragrafo si trovano elencati i seguenti esercizi che hanno come obiettivo la costruzione, mediante ''sezioni di Dedekind'', dei reali positivi a partire da razionali \vspace{5mm} positivi: \\
\textbf{42.6.} A \textbf{Dedekind cut} is a nonempty proper subset $A$ of $\Q_+^{^*}$ having \\
the following two properties:\\
$1^\circ$ For all $x,y\in\Q_+^{^*}$ if $x<y$ and if $y\in A$, then $x\in A$.\\
$2^\circ$ For all $x\in\Q_+^{^*}$ if $x\in A$, then there exists $y\in A$ such that $y>x$.\\
Let $\CC$ be the set of all Dedekind cuts. For each $r\in\Q_+^{^*}$ we define $L_r$ by
$$L_r=\{x\in\Q_+^{^*}: x<r\}$$
$(a)$ If $A\in\CC$, then every strictly positive rational not belonging to $A$ is greater than every number of $A$. $(b)$ If $r\in\Q_+^{^*}$, then $L_r\in\CC$. $(c)$ The set $\CC$ is a stable subset of $(\CB(\Q), +_{\CB}, \cdot_{\CB})$\footnote{Si tratta dell'insieme delle parti di $\Q$ con le operazioni tra ''sistemi'' indotte dalle operazioni di $\Q$.}. We shall denote by $+$ and $\cdot$ the compositions induced on $\CC$ by those of $\CB(\Q)$. Thus $\CC$ is a commutative semigroup for both addition and multiplication. $(d)$ The Dedekind cut $L_1$ is the multiplicative identity element of $\CC$. $(e)$ Multiplication on $\CC$ is distributive over addition.\\
\textbf{42.7.} Every element $A$ of $\CC$ is invertible for multiplication. [Let $B=\{y\in\Q_+^{^*}: y<u^{-1} \text{for some}\; u\notin A\}$. To show that if $z<1$, then $z\in AB$, let $a\in A$ and consider integral multiples of $\df{a}{n}$ where $z<\df{n}{n+1}$].\\
\textbf{42.8.} $(a)$ The ordering $\subseteq$ on $\CC$ is a total ordering compatible with both addition and multiplication. $(b)$ If $\CA$ is a nonempty subset of $\CC$ and if $D$ is an element of $\CC$ satisfying $X\subseteq D$, for all $X\in\CA$, then $\cup\CA\in\CC$ and $\cup\CA$ is the supremum of $\CA$ in $\CC$.\\
\textbf{42.9.} $(a)$ If $A, B\in\CC$, then $A\subset A+B$. [Consider fractions whose denominator $n$ satisfies $\df{1}{n}\in B$.] (b) If $A$ and $C$ are elements of $\CC$ and if $A\subset C$, then there exists $B\in\CC$ such that $A+B=C$. [Let $b=\{c'-c : c'\in C, c'>c, \text{and}\; c\notin A\}$. If $c\in C$ and $c\notin A$, to show that $c\in A+B$, consider fractions whose denominator $n$ satisfies $c+\df{1}{n}\in C$.] $(c)$ Every element of $\CC$ is cancellable for addition.\\
\textbf{42.10.} Let $(\R, +)$ be an inverse-completion of the semigroup $(\CC, +)$. There exists a composition $\cdot$ and a total ordering $\leq$ on $\R$ such that $(\R, +, \cdot, \leq)$ is a Dedekind ordered \vspace{0.5cm} field. \\
In questo lavoro, seguendo lo schema dato dagli esercizi, presentiamo una costruzione dei reali che, a nostro parere, non solo presenta parecchi vantaggi in termini di chiarezza ed essenzialit\a ma che permette di valorizzare pienamente il ruolo delle strutture $\Q_+^{^*}$ e $\R_+^{^*}$.
\section{\large Costruzione di $\R_+^{^*}$}
\begin{definition} Chiamiamo segmento iniziale\footnote{Preferiamo parlare di segmento iniziale perch\e questo sottinsieme \e la prima parte di una sezione in $\Q^*_+$ nell'accezione usuale del termine.} di $\Q^*_+$ un sottoinsieme $A$ proprio
e non vuoto di $\Q_+^{^*}$ con le seguenti propriet\a:
\begin{enumerate}
\item[$\rm{(1)}$]  Se $y\in A$  e $x<y$, allora $x\in A$, 
\item[$\rm{(2)}$]  Se $x\in A$, allora esiste $y\in A$ tale che $y>x$.
\label{segm}
\end{enumerate}
\end{definition}
L'insieme $A=\{x\in\Q_+^{^*} : x^2<2\}$ fornisce un 
esempio di segmento iniziale. Dimostriamo che $A$ verifica la $(2)$ della definizione 
\ref{segm}. 
Sia $x=\df{a}{b}\in A$ con $a, b$ coprimi.  
Se $\df{a}{b}<1$, possiamo assumere $y=1$. \\Quindi supponiamo $a\geq b$. Da $\displaystyle\frac{a^2}{b^2}<2$ segue \vspace{2mm} $\displaystyle\frac{a}{b}<2$. \\
Costruiamo una successione di frazioni nel seguente \vspace{2mm} modo: 
il primo termine della successione \e la mediante\footnote{$(\Q_+^{^*}, <)$ \e denso; se $\displaystyle\frac{a}{b}, \ \frac{c}{d}\in \Q_+^{^*}$ \  e \  
$\displaystyle\frac{a}{b} <\frac{c}{d}$ si ha $\displaystyle\frac{a}{b}< \frac{a+c}{b+d}< \frac{c}{d}$, chiamiamo \it{mediante} la frazione $\displaystyle\frac{a+c}{b+d}$.  \vspace{2mm}\\} 
$\displaystyle\frac{a+2}{b+1}$  tra la \vspace{3mm} frazione $\displaystyle\frac{a}{b}$ 
 e la frazione $\displaystyle\frac{2}{1}$, \\
il secondo termine \e la mediante  $\displaystyle\frac{2a+2}{2b+1}$ 
tra \vspace{3mm} $\displaystyle\frac{a}{b}$ e $\displaystyle\frac{a+2}{b+1}$, 
procedendo cos\io\!\!, \\
l'ennesimo termine della successione \vspace{3mm} sar\a $\displaystyle\frac{na+2}{nb+1}$.\\
Per \vspace{3mm} ogni $n\in\N$,  
$\displaystyle\frac{na+2}{nb+1}>\displaystyle\frac{a}{b}$.\\
Mostriamo che, per $n$ sufficientemente grande,  esistono termini della successione che 
al quadrato sono minori di $2$, cio\e che \vspace{3mm} appartengono ad $A$. 
Infatti $(\displaystyle\frac{na+2}{nb+1})^2<2$ se e solo se $(na+2)^2<2(nb+1)^2$ che \vspace{3mm} equivale a 
$n[n(2b^2-a^2)+4(b-a)]-2>0$.\vspace{3mm}
Se prendiamo 
$n>\displaystyle\frac{4(a-b)}{2b^2-a^2}$ si ottiene  $n(2b^2-a^2)+4(b-a)\geq 1$,  pertanto,
se \e anche  $n\geq 3$,
 la frazione corrispondente \e una di quelle cercate.
\begin{proposition}
Se $A$ \e un segmento iniziale di $\Q_+^{^*}$ e $B=\Q_+^{^*}\backslash A$, allora $A$ e $B$ formano una sezione 
in $\Q_+^{^*}$.
\label{sezionerepo}
\end{proposition}
$Dimostrazione.$  Dobbiamo solo verificare che se $z\in B$, allora
$x<z$ per ogni $x\in A$. Se vi fosse un $y\in A$ con $y>z$, allora,
per la $(1)$ della definizione \ref{segm}, l'elemento 
$z$ dovrebbe stare in $A$. \vspace{2mm} \bx
\par\noindent 
Chiameremo $A|B$ sezione associata al segmento iniziale $A$. Osserviamo \vspace{1mm} che:
\par\noindent
$\bullet$ Se la sezione associata ad un segmento iniziale di $\Q_+^{^*}$  \e di Dedekind, per la 
definizione \ref{segm}, l'elemento separatore della sezione \e necessariamente il minimo dell'insieme $B$. \\
$\bullet$ Per ogni $r\in \Q_+^{^*}$, sia $S_r=\{x\in \Q_+^{^*} : x<r \}$; essendo $\Q_+^{^*}$ denso,
$S_r$ \e un segmento iniziale. 
La sezione $A|B$ in $\Q_+^{^*}$, 
\begin{equation}
A=S_r  \ \ \text{e} \ \ B=\Q_+^{^*}\backslash S_r=\{x\in \Q_+^{^*} : x\geq r \},
\end{equation} 
associata al segmento iniziale $S_r$ 
 \e una sezione di Dedekind con elemento separatore $r$. \\
$\bullet$ Di conseguenza le sezioni di Dedekind in $\Q_+^{^*}$ associate a segmenti iniziali di $\Q_+^{^*}$ 
sono tutte e sole quelle del tipo \vspace{2mm} $S_r|(\Q_+^{^*}\backslash S_r)$.
\par\noindent
Chiamiamo insieme dei numeri reali positivi $\R_+^{^*}$ l'insieme di tutti i
segmenti iniziali di \vspace{4mm} $\Q_+^{^*}$. \\
\textbf{Operazioni} \vspace{3mm} \\
Usare i segmenti iniziali, anzich\e le sezioni in $\Q_+^{^*}$, permette di introdurre
con maggiore semplicit\a le operazioni in $\R_+^{^*}$.
Comunque presi $A, A'\in \R_+^{^*}$, definiamo:
\begin{equation}
A+A'=\{x+x' : x\in A, x'\in A'\}
\label{addrep}
\end{equation}
\begin{equation}
AA'=\{xx' : x\in A, x'\in A'\}
\label{moltre}
\end{equation}
\begin{proposition}
$\R_+^{^*}$ \e chiuso rispetto alle operazioni definite in ${(\ref{addrep})}$ e ${(\ref{moltre})}$.
\end{proposition}
$Dimostrazione.$ Dobbiamo provare che $A+A'$ e $AA'$ sono segmenti
iniziali di $\Q_+^{^*}$.
Se $z\notin A$ e $z'\notin A'$, allora $z+z'\notin A+A'$ e $zz'\notin AA'$, infatti, per la definizione di segmento iniziale,
per ogni $x\in A$ e per ogni $x'\in A'$, si ha $x<z$ e $x'<z'$  da
cui segue $x+x'<z+z'$ e $xx'<zz'$; cos\io $A+A'$ e $AA'$ sono sottoinsiemi
propri e non vuoti di $\Q_+^{^*}$.\\
Proviamo la $(1)$ della definizione \ref{segm}. Sia $y+y'\in A+A'$ e \vspace{1mm} sia 
$w<y+y'$; consideriamo $t=\displaystyle\frac{wy}{y+y'}$ e $t'=\displaystyle\frac{wy'}{y+y'}$, si \vspace{1mm} ha $t+t'=w$.\\ 
L'asserto segue osservando che $t<y$ e $t'<y'$ implicano, essendo
$A$ e $A'$ segmenti iniziali, $t\in A$ e $t'\in A'$.\\
Sia ora $yy'\in AA'$ e $w<yy'$, allora $\displaystyle\frac{w}{y'}<y$. Essendo $\Q_+^{^*}$ denso,
 \vspace{1mm} esiste $v\in \Q_+^{^*}$ tale che $\displaystyle\frac{w}{y'}<v<y$ quindi $v\in A$.
Posto $v'=\displaystyle\frac{w}{v}$,  si ha $v'<y'$, cio\e $v'\in A'$ e $w=vv'$.\\
La $(2)$, della definizione \ref{segm}, per $A+A'$ e $AA'$ segue 
dal fatto che tale propriet\a vale per $A$ e per $A'$. \vspace{2mm} \bx
\par\noindent
Le propriet\a associativa e commutativa dell'addizione in $\Q_+^{^*}$ ci permettono di verificare che:
\begin {proposition}
L'operazione di addizione in $\, \R_+^{^*} \,$ gode delle propriet\a commutativa e associativa.
\label{comassrepo}
\end{proposition}
\begin{lemma}
Fissati $A\in \R_+^{^*}$ e $n\in \N$, esistono $x\in A$ e $y\notin A$ tali \vspace{1mm} che 
$x+ \displaystyle\frac{1}{n}>y$.
\label{lemmare}
\end{lemma}
$Dimostrazione.$ Siano $x_0\in A$ e $y_0\notin A$ cosicch\ea\ $y_0=x_0+t$,
per un \vspace{1mm} opportuno $t\in \Q_+^{^*}$. Scegliamo un naturale\footnote{La
propriet\a di Archimede nell'insieme dei numeri naturali dice che:
\textit{Quali che siano $a$ e $b$ in $\N$, esiste un naturale $n$ tale che $na>b$.}
In $\Q_+^{^*}$ la propriet\a pu\o essere cos\io riformulata 
$\forall r\in\Q^{^+} \ \text{esiste} \ n\in\N \ \text{tale che} \ n>r$:
Se $r=\df{a}{b}$, essendo $b\geq 1$, si ha $ab+b>a$ ovvero $a+1>\df{a}{b}$.} $k>nt$ e
poniamo $v= \displaystyle\frac{t}{k}$.\\
Consideriamo i numeri 
$$x_0, \ x_0+v, \ x_0+2v,\  ..., \  x_0+(k-1)v, \ x_0+kv=x_0+t=y_0.$$
\E evidente, essendo $x_0\in A$ e $y_0\notin A$, che vi \e un $j$,
$1\leq j\leq k$, per cui $x_0+(j-1)v\in A$ e $x_0+jv\notin A$. Posto
$x=x_0+(j-1)v$ e $y=x_0+jv$, si ottiene:
$$y=x_0+jv=x_0+(j-1)v+v=x+v=x+ \displaystyle\frac{t}{k}<x+ \displaystyle\frac{t}{nt}=x+ \displaystyle\frac{1}{n}.$$ \bx

\begin{corollary}\label{corolre}
Fissati $A\in \R_+^{^*}$ e $m\in \N$, esistono $x\in A$ e $y\notin A$ tali \vspace{1mm} che
$\displaystyle\frac{x}{y}+ \displaystyle\frac{1}{m}>1$.
\end{corollary}
$Dimostrazione.$ Sia $x_1\in A$ e scegliamo un $h\in \N$ tale \vspace{2mm} che $h> \displaystyle\frac{m}{x_1}$.\\
Per il lemma precedente esistono $x\in A$ e $y\notin A$ tali \vspace{2mm} che $y<x+ \displaystyle\frac{1}{h}$,
da cui $1< \displaystyle\frac{x}{y}+ \displaystyle\frac{1}{hy}$. 
Avremo l'asserto se mostriamo \vspace{2mm} che 
$\displaystyle\frac{1}{hy}< \displaystyle\frac{1}{m}$.\\
Poich\ea\ $x_1\in A$ e $y\notin A$, per definizione di segmento \vspace{2mm} iniziale
$ \displaystyle\frac{x_1}{y}<1$, allora, essendo $ \displaystyle\frac{1}{h}< \displaystyle\frac{x_1}{m}$,
 si ha $\df{1}{hy}<\df{x_1}{my}<\df{1}{m}$. \vspace{5mm} \bx
\par\noindent
Siamo ora in grado di provare che:
\begin{proposition}\label{gruppo} 
$\R_+^{^*}$ \e, rispetto alla moltiplicazione, un gruppo abeliano.
\end{proposition}
$Dimostrazione.$ Le propriet\a commutativa e associativa seguono, anche qui, dalle analoghe propriet\a in $\Q_+^{^*}$.\\  
Il segmento iniziale $S_1=\{x\in \Q_+^{^*} : x<1\}$ di $Q_+^{^*}$ \e l'elemento neutro,
proviamo che, per ogni $A\in \R_+^{^*}$, $S_1A=A$.\\ 
Prendiamo un qualunque elemento $xx'\in S_1A$. Da $x<1$ segue
$xx'<x'$ e quindi $xx'\in A$.
Dimostriamo ora che $A\subset S_1A$. Dalla $(2)$ della definizione \ref{segm} segue che,
per ogni $y\in A$, esiste $x\in A$ tale che $y<x$. Allora \vspace{1.5mm} possiamo scrivere
$y=(\displaystyle\frac{y}{x})x$ con $\displaystyle\frac{y}{x}<1$.\\
Proviamo che esiste l'inverso di ogni
elemento di  $\R_+^{^*}$. 
Preso comunque $A\in \R_+^{^*}$,
poniamo $\widehat{A}=\{x\in \Q_+^{^*} : x< \displaystyle\frac{1}{y}, \text{per qualche} \; y\notin A\}$.\\
Dimostriamo che $\widehat{A}$ \e un segmento iniziale di $\Q_+^{^*}$.\\
Essendo il complementare di $A$ in $\Q_+^{^*}$ non vuoto, anche $\widehat{A}$ \e non \vspace{1mm} vuoto.\\
Preso comunque
$y\in A$ , ogni elemento $z> \displaystyle\frac{1}{y}$ non \e in $\widehat{A}$, essendo
$\df{1}{y}>\df{1}{t}$ per ogni $t\notin A$, quindi $\widehat{A}$ \e un sottoinsieme proprio di $\Q_+^{^*}$.\\
Dimostriamo ora le due propriet\a della definizione \ref{segm}. 
La prima \e ovvia per come abbiamo definito $\widehat{A}$.\\
 Sia ora $x\in\widehat{A}$,
allora $x<\displaystyle\frac{1}{z}$ per qualche $z\notin A$. Essendo $\Q_+^{^*}$ denso,
esiste $t\in\Q_+^{^*}$ tale che $x<t<\df{1}{z}$; l'elemento $t$ \e quello richiesto dalla \vspace{1.5mm} definizione.\\
Proviamo adesso che $A\widehat{A}=S_1$. Se $z\in A\widehat{A}$, allora 
$z=x\hat{x}$ con $x\in A$ e \vspace{1mm} $\hat{x}\in\widehat{A}$.
Quindi $z=x\hat{x}< \displaystyle\frac{x}{y}$ per qualche $y\not\in A$; pertanto,
essendo $x\in A$ e $y\notin A$, $z<\displaystyle\frac{x}{y}<1$ cio\e $z\in S_1$.\\
Viceversa sia $w\in S_1$, e sia $m\in \N$ tale che\footnote{Essendo $\Q_+^{^*}$ denso, esiste $\frac{n}{m}$ con $w<\frac{n}{m}<1$. Poich\e $n<m$, si ha 
$n\leq m-1$ e $\frac{n}{m}\leq \frac{m-1}{m}$.}
 $w\leq \displaystyle\frac{m-1}{m}<1$. \\
Per il corollario precedente, esistono $x\in A$ e $y\notin A$ tali \vspace{1mm} che
$\displaystyle\frac{x}{y}>\displaystyle\frac{m-1}{m}$ da cui 
$\displaystyle\frac{1}{y}>\displaystyle\frac{m-1}{mx}$. Poich\ea\ $\Q_+^{^*}$ \e denso, esiste
$v\in \Q_+^{^*}$\vspace{1mm} con $\displaystyle\frac{m-1}{mx}<v<\displaystyle\frac{1}{y}$, ovvero
$\displaystyle\frac{m-1}{m}<vx<\displaystyle\frac{x}{y}$ e quindi $w<vx$ con $x\in A, v\in\widehat{A}$ . \vspace{5mm} \bx
\par\noindent
Infine dimostriamo che:
\begin{proposition}\label{distrepo}
In $\R_+^{^*}$ vale la propriet\a distributiva della moltiplicazione rispetto all'addizione:
 $A(A'+A'')=AA'+AA''$, per ogni $A, A', A''\in \R_+^{^*}$.
\end{proposition}
$Dimostrazione.$ L'inclusione $A(A'+A'')\subseteq AA'+AA''$ segue dalla 
propriet\a distributiva in $\Q_+^{^*}$. \\
Preso comunque $xx'+yx''\in AA'+AA''$, possiamo supporre, senza ledere
le generalit\ag, che $x\geq y$. Da
$xx'+yx''\leq xx'+xx''=x(x'+x'')$ segue, per la definizione
\ref{segm} applicata al segmento iniziale $A(A'+A'')$, che $xx'+yx''\in A(A'+A'')$.
Ci\o prova l'inclusione $AA'+AA''\subseteq A(A'+A'')$. \vspace{5mm} \bx
\par\noindent
\textbf{Alcune propriet\a di $\R_+^{^*}$}\begin{proposition}
L'inclusione in $\R_+^{^*}$ \e una relazione di ordine totale compatibile con l'addizione 
e la moltiplicazione.
\label{ordrepo}
\end{proposition}
$Dimostrazione.$ La relazione $\subseteq$ \e una relazione d'ordine nell'insieme delle
parti di $\Q_+^{^*}$. La relazione da essa indotta su un sottoinsieme di parti di $\Q_+^{^*}$,
nel nostro caso in $\R_+^{^*}$,  \e ancora d'ordine, ovvero gode delle propriet\a riflessiva, antisimmetrica e transitiva.\\
Dimostriamo che questa relazione d'ordine \e totale. Supponiamo 
$A\nsubseteq A'$, allora esiste $x\in A$ tale che $x\notin A'$. 
Per definizione di segmento iniziale, per ogni $y\in A'$ si ha $y<x$ che implica 
$y\in A$ e quindi $A'\subset A$.\\
\'E facile verificare la compatibilit\a di questa relazione con le
operazioni in $\R_+^{^*}$, cio\e\!\!:
\begin{equation}
\text{Se} \; A\subseteq A' \; \text{allora} \; A+A''\subseteq A'+A'', \; \forall A''\in\R_+^{^*}
\label{comparead}
\end{equation}
\begin{equation}
\text{Se} \; A\subseteq A' \; \text{allora} \; AA''\subseteq A'A'', \; \forall A''\in\R_+^{^*}
\label{comparemolt}
\end{equation}
\bx
\begin{proposition}
$\R_+^{^*}$ \e un insieme continuo.
\label{continuo}
\end{proposition}
$Dimostrazione.$ Nella proposizione precedente abbiamo visto
che $(\R_+^{^*}, \subseteq )$ \e un insieme totalmente ordinato, dobbiamo quindi provare che \e denso ed ha la propriet\a di Dedekind.\\
Siano $A, A'\in \R_+^{^*}$ con $A\subseteq A'$, $A\neq A'$ e sia $x'$ un elemento in $A'$
 non appartenente ad $A$.\\
Essendo $A'$ un segmento iniziale, esiste $r\in A'$ tale che $r>x'$, da cui 
$A\subseteq S_r$ e $A\neq S_r$, perch\ea\ $a>r$, per qualche $a\in A$, implicherebbe $a>x'$ ovvero $x'\in A$.  
Inoltre, essendo $r\in A'$, $S_r\subseteq A'$, e poich\e esiste $y\in A'$ con $y>r$,  $S_r\neq A'$.\\
Sia ora $\CA$ un sottoinsieme non vuoto di $\R_+^{^*}$ limitato superiormente; mostriamo che
$\CA$ ammette estremo superiore.\\
Poniamo 
$$\widehat{A}=\bigcup_{A\in\CA} A.$$
Sia $\tilde{A}\in\R_+^{^*}$ una limitazione superiore per $\CA$. Da $A\subseteq \tilde{A}$, per ogni $A\in \CA$,
segue che $\widehat{A}\subseteq\tilde{A}$, cosicch\ea\ $\widehat{A}$ \e un sottinsieme proprio (e non vuoto) 
di $\Q_+^{^*}$. Si verifica facilmente che $\widehat{A}$ \e un segmento iniziale di $\Q_+^{^*}$. \\
$\widehat{A}$ \e l'estremo superiore cercato, essendo $A\subseteq\widehat{A}\subseteq\tilde{A}$ per ogni 
$A\in \CA$ e per ogni limitazione superiore $\tilde{A}$ di $\CA$.\bx

\vspace{0.2cm}
\par\noindent
Mostriamo ora che in $\R_+^{^*}$ vale la legge di cancellazione per l'addizione. 
\begin{lemma} 
Per ogni $A, A'\in \R_+^{^*}$, si ha  $A\subsetneq A+A'$. 
\label{lemmare2}
\end{lemma}
$Dimostrazione.$ Preso comunque $x\in A$, si ha $x<x+x'$ per ogni $x'\in A'$, quindi $x\in A+A'$. \\
Mostriamo che l'inclusione \e propria. Sia $z\in A'$ e scegliamo \vspace{1mm} un $n\in \N$ tale che 
$\df{1}{n}<z$.
Per il lemma \ref{lemmare}, esistono $x\in A$ e $y\notin A$ tali \vspace{1mm} che 
$x+\df{1}{n}>y$. Pertanto $x+z>y$ e quindi $x+z\notin A$.\bx

\begin{proposition}
Per ogni $A, A'\in \R_+^{^*}$ con $A\subsetneq A'$, esiste $A''\in \R_+^{^*}$  tale che $A'= A+A''$. 
\label{complerepo}
\end{proposition}
$Dimostrazione.$ Poniamo
 $$A''=\{z\in \Q_+^{^*} : z+y=x \; \text{con} \; x\in A', y\notin A, x>y\};$$
proviamo che questo sottoinsieme di $\Q_+^{^*}$ \e un segmento iniziale.\\
$A''$ non \e vuoto. Infatti, 
essendo $A\neq A'$, esiste $y\in A'$ con $y\notin A$. Poich\e 
esiste $x\in A'$ con $x>y$, l'elemento $z$ tale che 
$z+y=x$ \e un elemento di $A''$.\\
Mostriamo che $A''\neq \Q_+^{^*}$.
Sia $w\notin A'$. Poich\e, per ogni $y\in A'$ $y\notin A$, $w+y>w$ si ha
$w+y\notin A'$ che implica $w\notin A''$.\\
Restano da provare la $(1)$ e $(2)$ della definizione \ref{segm}.\\
Sia $z\in A''$, allora possiamo trovare $x\in A'$ e $y\notin A$ con $x>y$ tali che 
$z+y=x$. Se $v<z$,  da $x_1=v+y<z+y=x$, segue $x_1\in A'$ e quindi $v\in A''$.\\
Infine preso comunque $z\in A''$, allora $z+y=x$ con $x\in A'$, $y\notin A$ e 
$x>y$. Essendo $A'$ un segmento iniziale, esiste $x'\in A'$ tale che $x'>x$
quindi con $x'>y$. Consideriamo $t\in A''$ con $t+y=x'$; da 
$t+y=x'>x=z+y$ si ottiene, in $\Q_+^{^*}$, $t>z$.\\
Proviamo infine che $A'=A+A''$.\\
Cominciamo con l'inclusione $A+A''\subseteq A'$.
Prendiamo $t+z\in A+A''$, allora esistono $x\in A'$ e $y\notin A$ 
con $x>y$ tali che $z+y=x$. Essendo $t<y$, si ha 
$t+z<z+y=x$ che implica $t+z\in A'$.\\
Viceversa, $A'\subseteq A+A''$.
Sia $y\in A'$. Se $y\in A$, per il lemma \ref{lemmare2}, si ha $y\in A\subsetneq A+A''$.
Quindi possiamo supporre che $y\notin A$. Giacch\ea\ 
esiste $x'\in A'$ tale che $x'>y$, possiamo considerare
l'elemento $z\in A''$ tale che $x'=y+z$.\\
Sia $n\in \N$ tale che $\df{1}{n}<z$; per il lemma \ref{lemmare}, esistono $t\in A$ e
$v\notin A$ tale che $t+\df{1}{n}>v$. Essendo $y>t$, si ha
$x'=y+z>y+\df{1}{n}>t+\df{1}{n}>v$ e quindi si pu\o considerare
l'elemento $w\in A''$ tale che $w+v=x'$.\\
Osserviamo che $t+\df{1}{n}+w>v+w=x'=y+z>y+\df{1}{n}$, da cui si ottiene $t+w>y$ con $t\in A$ e $w\in A''$ e possiamo 
concludere che $y\in A+A''$.\bx

\begin{corollary}
In $\R_+^{^*}$ vale la legge di cancellazione
\label{canrepo}
\end{corollary}
$Dimostrazione.$ Dobbiamo mostrare che
$$A+\bar{A}=A'+\bar{A} \; \;  \text{implica} \; \; A=A'.$$
\par\noindent
Supponiamo per assurdo che sia  $A\neq A'$; essendo $\R_+^{^*}$ un insieme 
totalmente ordinato, possiamo assumere $A\subseteq A'$. Per la proposizione
precedente, esiste $A''$ tale che $A'=A+A''$.
Allora si ha $A+\bar{A}=A'+\bar{A}=A+\bar{A}+A''$ che contraddice il lemma
\ref{lemmare2}.\bx

\begin{corollary}
 $A\subsetneq A'$, implica $A+\bar{A}\subsetneq A'+\bar{A}$.
\end{corollary}

\par\noindent
Concludiamo il paragrafo provando che esiste un monomorfismo da $\Q_+^{^*}$ in $\R_+^{^*}$.
Consideriamo l'applicazione 
$\varphi : \Q_+^{^*}\longrightarrow \R_+^{^*} $ definita da $\varphi(r)=S_r$.
Siano $r, r'\in \Q_+^{^*}$ con $r < r'$. Essendo $\Q_+^{^*}$ denso, esiste un razionale positivo $t$ tale che $r<t<r'$ da cui
$S_r\subsetneq S_{r'}$, quindi $\varphi$ \e un'applicazione iniettiva e isotona (conserva l'ordine).\\
Dimostriamo che $\varphi$ \e un omonorfismo: per ogni $r, r'\in \Q_+^{^*}$ deve essere 
$$\varphi(r+r')=\varphi(r)+\varphi(r')\; \; \text{e} \; \; \varphi(rr')=\varphi(r)\varphi(r')$$
\par\noindent
ovvero 
$$S_{r+r'}=S_r+S_{r'} \; \; \text{e} \; \;  S_{rr'}=S_rS_{r'}.$$
\par\noindent
Presi $x\in S_r, x'\in S_{r'}$, si ha $x+x'<r+r'$ e quindi $S_r+S_{r'}\subseteq S_{r+r'}$; inoltre
\vspace{0.2cm}
\par\noindent
potendo scrivere ogni 
$ x\in S_{r+r'}$ nella forma $x=y+y'$ con $y=\displaystyle\frac{xr}{r+r'}\in S_r$, 
e  $y'=\displaystyle\frac{xr'}{r+r'}\in S_{r'}$, \e vera anche l'inclusione $S_{r+r'}\subseteq S_r+S_{r'}$
\vspace{0.2cm}
\par\noindent
Come sopra, da $xx'<rr'$, segue $S_rS_{r'}\subseteq S_{rr'}$; preso comunque un 
\vspace{1mm} elemento $x\in S_{rr'}$, 
essendo $\displaystyle\frac{x}{r'}<r$ ed essendo $\Q_+^{^*}$ denso, esiste 
$v\in\Q_+^{^*}$ tale 
\vspace{0.2cm}
\par\noindent
che $\displaystyle\frac{x}{r'}<v<r$. 
Posto $v'=\displaystyle\frac{x}{v}$, si ha $v\in S_r$, $v'\in S_{r'}$ 
e $x=vv'$.
\section{\large Costruzione di $\R$}
Introduciamo nel prodotto cartesiano $\mathbb{R}^{^*}_+\times\mathbb{R}^{^*}_+$ la seguente relazione binaria:
\begin{equation}\label{eqreal}
(A, A')\approx (B, B')  \stackrel{def}{\Longleftrightarrow} A+B'=A'+B
\end{equation}
Usando le propriet\a di cui gode l'addizione in $\R_+^{^*}$, si dimostra che la 
relazione data \e di equivalenza.\\ 
Chiamiamo insieme dei numeri reali $\mathbb{R}$ l'insieme quoziente 
$(\mathbb{R}^{^*}_+\!\times\mathbb{R}^{^*}_+)/\approx$\\
Presi gli elementi $[(A, B)], [(A', B')]\in \mathbb{R}$, poniamo
\begin{equation}\label{addre}
[(A, B)] + [(A', B')]  \stackrel{def}{=}[(A+A', B+B')]
\end{equation}
\begin{equation}\label{molre}
[(A, B)]\cdot  [(A', B')]\stackrel{def}{=}[(AA'+BB', AB'+BA')]
\end{equation}
\par\noindent
Tali operazioni sono ben definite, nel senso che non dipendono dalla scelta dei rappresentanti. Se $[(A, B)]=[(A', B')]$ \ e \ $[(C, D)]=[(C', D')]$, \ abbiamo le seguenti uguaglianze 
\ $(A+B')C=(B+A')C$, \ $(B+A')D=(A+B')D$, \ $A'(C+D')=A'(D+C')$ \ e
 $B'(D+C')=B'(C+D')$, \ da cui, applicando la propriet\a distributiva in $\mathbb{R}^{^*}_+$ e sommando membro a membro, si ottiene  
$$AC+B'C+BD+A'D+A'C+A'D'+B'D+B'C'=$$
$$=BC+A'C+AD+B'D+A'D+A'C'+B'C+B'D'.$$ 
Per la legge di cancellazione in $\mathbb{R}^{^*}_+$, 
 $$AC+BD+A'D'+B'C'=AD+BC+A'C'+B'D',$$ 
da cui, per la {(\ref{eqreal})} e la {(\ref{molre})}, si ha
\ $[(A, B)] \cdot [(C, D)]=[(A', B')] \cdot [(C', D')]$. Analogamente si procede per l'addizione.
\begin{proposition}
$(\R, +, \cdot)$ è un dominio d'integrit\a con unit\a\!\!.
\end{proposition}
$Dimostrazione$. \'E facile verificare che l'addizione e la moltiplicazione godono delle propriet\a commutativa e associativa e che vale la propriet\a distributiva della moltiplicazione rispetto all'addizione. \\
Osserviamo che  tutte le coppie con le due componenti uguali tra loro sono
in relazione e scegliamo come rappresentante di questa classe la coppia $(S_1, S_1)$.
Per ogni $[(A, B)]\in \R$,
$$[(A, B)]+[(S_1, S_1)]=[(A+S_1, B+S_1)]=[(A,B)]$$
\par\noindent 
da cui segue che la classe \ $[(S_1, S_1)]$ \ \`e l'elemento neutro rispetto all'addizione.\\
Inoltre l'elemento \ $[(B, A)]\in \R$ \ \`e 
tale che
$$[(A, B)]+[(B, A)]= [(A+B, B+A)]=[(S_1, S_1)]$$
ovvero $[(B, A)]$ \`e l'opposto di $[(A, B)]$. \\
Infine si ha: 
$$[(A, B)] \cdot [(C+S_1, C)]=[(AC+A+BC, AC+BC+B)]=[(A, B)],$$
 \par\noindent
e poich\ea, per ogni \ $C, \ D\in\mathbb{R}^{^*}_+$,  $(D+S_1, D)\approx (C+S_1, C)$, 
possiamo scegliere come 
rappresentante dell'elemento neutro rispetto alla moltiplicazione la coppia 
$(S_1+S_1, S_1)$.\\
Proviamo che in $\R$ vale la legge dell'annulamento del prodotto.\\
Sia\ $[(A, B)]\not =[(S_1, S_1)]$ \ e \ $[(A, B)] \cdot [(C, D)] = [(S_1, S_1)]$ \ cio\e, per la {(\ref{eqreal})}, \\
$A\not= B$ e $AC+BD = AD+BC$. Se $A\supset B$, per un opportuno reale positivo $H$ si ha $A=B+H$ (prop. \ref{complerepo}) \ da cui, sostituendo,
$BC+HC+BD = BD+HD+BC$ \ e quindi \ $C=D$ (cor. \ref{canrepo} e prop \ref{gruppo}). 
Si procede allo stesso modo se  $A\subset B$.\bx
\vspace{2mm}
\par\noindent
Prendiamo un elemento $[(A, B)]\in \R$. Se $A\supset B$
esiste un $C\in \R_+^{^*}$ (prop. \ref{complerepo}) tale che $A=B+C$; allora come rappresentante dell'elemento $[(A, B)]$
possiamo prendere la coppia $(S_1+C, S_1)$. Analogamente nel caso $A\subset B$, si pu\o 
scegliere come rappresentante della classe $[(A, B)]$ la coppia $(S_1, S_1+D)$, per un opportuno  
$D\in \R_+^{^*}$. Quindi possiamo scrivere:
\begin{equation} \label{int}
[(A, B)] = \left\{  \begin{array}{ccc}
[(S_1+C, S_1)] &  \textrm{se $A\supset B $} \cr
[(S_1, S_1)] &  \textrm{se $A=B$}  \cr
[(S_1, S_1+D)] &  \textrm{se $A\subset B$}
\end{array} \right.
\end{equation}
\par\noindent
\begin{proposition}
$\R$ \e un campo.
\end{proposition}
$Dimostrazione$. Preso il segmento iniziale $A\in \R_+^{^*}$ e indicato con $A'$ il suo inverso (prop. \ref{gruppo}),
\e facile verificare che l'inverso in $\R$ di un elemento del tipo \ $[(S_1+A, S_1)]$ \ (rispett. $[(S_1, S_1+A)]$) \
\e \ $[(S_1+A', S_1)]$ \ (rispett. \\
$[(S_1, S_1+A')]$).\bx
\par\noindent
Si prova che in $\R$ la seguente relazione binaria:
\begin{equation}\label{ordre}
[(A, B)] < [(A', B')]  \stackrel{def}{\Longleftrightarrow} A+B'\subset B+A'
\end{equation}
\par\noindent
\e una relazione d'ordine stretto tricotomica e che la relazione d'ordine $\leq$ ad essa associata
\begin{equation}\label{ordtotre}
[(A, B)] \leq [(C, D)]  \stackrel{def}{\Longleftrightarrow} [(A, B)] < [(C, D)]\; \text{oppure}\;
 [(A, B)] = [(C, D)]
\end{equation}
\e d'ordine totale.\\
Chiameremo \it{reali positivi} i numeri $[(A, B)]> [(S_1, S_1)]$ ovvero le classi del tipo $[(C+S_1, S_1)]$, 
\it{reali negativi} i numeri $[(A, B)]< [(S_1, S_1)]$ cio\e le classi $[(S_1, D+S_1)]$.
\vspace{0.5cm}
\par\noindent
Consideriamo l'applicazione iniettiva
$f :\mathbb{R_+^{^*}}\longrightarrow \mathbb{R}$ definita nel seguente modo:
$$\forall A\in\mathbb{R}^{^+} \;\; f(A) = [(A+S_1, S_1)]$$
Mostriamo che si tratta di un omomorfismo.
Comunque presi $A, B\in \mathbb{R_+^{^*}}$ si ha:
$$f(A+B) = [(A+B+S_1, S_1)] = [(A+B+S_1+S_1, S_1+S_1)]=$$
$$= [(A+S_1, S_1)] + [(B+S_1, S_1)] = f(A) + f(B)$$
\par\noindent
e ancora,
$$f(AB) = [(AB+S_1, S_1)] = [(AB+A+B+S_1+S_1, A+B+S_1+S_1)]=$$
$$=[(A+S_1, S_1)]\cdot [(B+S_1, S_1)] = f(A)\cdot f(B)$$
\par\noindent
Osserviamo che l'applicazione $f$ \e anche isotona.\\
Essendo la sottostruttura di $\R$ formata dai reali positivi isomorfa ad $\R_+^{^*}$, possiamo identificare la classe $[(A+S_1, S_1)]$ con  $A$, allora la classe opposta 
$[(S_1, A+S_1)]$ sar\a indicata con $-A$ e infine 
per la classe $[(S_1, S_1)]$ useremo il simbolo $0$ (lo zero).\\
Con le notazioni introdotte, ogni numero reale pu\o essere scritto come differenza di reali positivi\footnote{Qui il problema di immersione \e il seguente:
Dato $H(+)$ un mezzo gruppo abeliano (un insieme in cui \e definita un'operazione che gode delle propriet\a commutativa e associativa e in cui vale la legge di cancellazione), costruire un gruppo abeliano $K(+)$, contenente un mezzo gruppo $H'(+)$ isomorfo ad $H(+)$, e tale che ogni elemento $\xi$ di $K$ possa esprimersi come \textit{differenza}, $\xi=a'-b'$, di due elementi di $H'$. Quest'ultima propriet\a permette di dimostrare che $K(+)$ \e il pi\u piccolo gruppo abeliano che contiene $H(+)$.  Il gruppo $K$ si chiama \textit{gruppo delle differenze} di $H$. Cfr. [4]}, infatti:
$$[(A, B)] = [(A+S_1+S_1, B+S_1+S_1)] = [(A+S_1, S_1)] + [(S_1, B+S_1)] = A-B$$
\par\noindent
Mostriamo ora che:
Se $(G, +)$ \`e un gruppo abeliano che contiene una sottostruttura isomorfa a 
$(\mathbb{R_+^{^*}}, +)$ allora $G$ contiene anche un sottogruppo  isomorfo a 
$(\mathbb{R}, +)$.\\
Supponiamo quindi che esistano un gruppo abeliano $G$ e un monomorfismo 
$g:\mathbb{R_+^{^*}}\longrightarrow G$ e definiamo l'applicazione $g^*:\mathbb{R}\longrightarrow G$
ponendo\footnote{Con $-g(B)$ indichiamo l'opposto di $g(B)$ in $G$.}
$$g^*([(A, B)])=g(A)-g(B).$$
L'applicazione $g^*$ \e ben definita ed iniettiva. 
Inoltre
$$g^*([(A, B)]+[(C, D)])=g^*([(A+C, B+D)])=g(A+C)-g(B+D)=$$
$$=g(A)-g(B)+g(C)-g(D)=g^*([(A, B)])+g^*([(C, D)]),$$
cosicch\ea\ $g^*$ \e un monomorfismo.\\
L'immersione di $\Q$ in $\R$ \e data dal monomomorfismo
$g :\mathbb{Q}\longrightarrow \mathbb{R}$ definito da 
$\;\; g[(r, r')] = [(S_r, S_{r'})]$.
\vspace{0.2cm}
\par\noindent
Concludiamo mostrando che:
\begin{proposition}
$\R$ \e continuo.
\end{proposition}
$Dimostrazione$. E' sufficiente verificare che ogni sezione di $\R$  \e di Dedekind.\\
Sia $(\CA, \CB)$ una sezione di $\R$. Se in $\CA$ ci sono elementi positivi
l'insieme delle controimmagini mediante l'applicazione $f :\mathbb{R}^{^+}\longrightarrow \mathbb{R}$ (l'immersione da $\R_+^{^*}$ in $\R$ data prima), 
 di tali elementi \e limitato superioremente in $\R_+^{^*}$ dalla controimmagine di qualche elemento di $\CB$.
Poich\e $\R_+^{^*}$ ha la propriet\a di Dedekind (prop. \ref{continuo}), questo insieme ammette estremo superiore la cui
immagine \e il minimo di $\CB$, ossia l'elemento separatore della classe.\\
Supponiamo ora che in $\CB$ ci siano elementi negativi.
Consideriamo la sezione $(\CB', \CA')$ di $\R$, dove $\CB'$, $\CA'$ contengono, rispettivamente,
tutti gli elementi opposti agli elementi di $\CB$ e $\CA$. Per quanto detto esiste
il minimo di $\CA'$ elemento separatore della sezione $(\CB', \CA')$, l'opposto di questo elemento sar\a il massimo di
$\CA$ e l'elemento separatore della sezione $(\CA, \CB)$.\\
Supponiamo infine che in $\CA$ non ci siano positivi e in $\CB$ non ci siano negativi. 
L'elemento separatore della sezione \e in questo caso lo zero, infatti poich\ea\ gli insiemi 
$\CA$ e $\CB$ formano, per definizione, una partizione di $\R$, lo zero deve appartenere ad uno di \vspace{16mm} essi. \bx
\newpage\noindent
{\large \textbf{Bibliografia}} \vspace{6mm} \\
\textbf{[1]} C.B. Boyer, \textit{Storia della Matematica}, Arnoldo Mondadori, Milano, 1982. \\
\textbf{[2]} U. Bottazzini, \textit{Il Flauto di Hilbert}, UTET, Torino, 1996. \\
\textbf{[3]} M.R. Enea, D. Saeli, \textit{Sistemi Numerici, un'introduzione}, Aracne, \vspace{-1.3mm} Roma, \hspace*{5mm} 2009. \\
\textbf{[4]} L. Lombardo-Radice, \textit{Istituzioni di algebra astratta}, Feltrinelli, \vspace{-1.3mm} Milano, \hspace*{5mm} 1967. \\
\textbf{[5]} S. Warner, \textit{Modern algebra}, Dover, New York, 1990.
\end{document}